\documentclass[11pt]{article}
\usepackage{amssymb,amsfonts,amsmath,latexsym,epsf,tikz,url}

\newtheorem{theorem}{Theorem}[section]

\newtheorem{observation}[theorem]{Observation}

\newtheorem{problem}[theorem]{Problem}
\newtheorem{corollary}[theorem]{Corollary}

\newtheorem{remark}[theorem]{Remark}

\newcommand{\proof}{\noindent{\bf Proof. }}
\newcommand{\qed}{\hfill $\square$\medskip}

\begin{document}

\title{The distinguishing number (index)  and the  domination number  of a graph}

\author{
Saeid Alikhani  
\and
Samaneh Soltani $^{}$\footnote{Corresponding author}
}

\date{\today}

\maketitle

\begin{center}
Department of Mathematics, Yazd University, 89195-741, Yazd, Iran\\
{\tt alikhani@yazd.ac.ir, s.soltani1979@gmail.com}
\end{center}


\begin{abstract}
The distinguishing number (index) $D(G)$ ($D'(G)$) of a graph $G$ is the least integer $d$
such that $G$ has an vertex labeling (edge labeling)  with $d$ labels  that is preserved only by a trivial
automorphism. A set $S$ of vertices in $G$ is a dominating set of $G$ if every vertex of $V(G)\setminus S$ is adjacent to some vertex in $S$. The minimum
cardinality of a dominating set of $G$ is the domination number of $G$ and denoted by $\gamma (G)$. 
In this paper, we obtain some upper bounds for the distinguishing number and the distinguishing index of a graph based on its domination number.
\end{abstract}


\medskip
\noindent{\bf AMS Subj.\ Class.}: 05C25, 05C69

\section{Introduction and definitions}

 Domination in graphs is very well studied in graph
theory. The literature on this subject has been surveyed and detailed in the two books by Haynes, Hedetniemi, and
Slater \cite{Haynes1,Haynes2}.
For notation and graph theory terminology we in general follow \cite{Haynes1}. Specifically, let $G =(V,E)$ be a graph with
vertex set $V$ of order $n= |V|$ and edge set $E$ of size $m = |E|$, and let $v$ be a vertex in $V$. The \textit{open neighborhood} 
of $v$ is the set $N(v)=\{u \in V |~ uv \in E\}$ and the \textit{closed neighbourhood} of $v$ is $N[v]=\{v\}\cup N(v)$. For a set $S$
of vertices, the open neighbourhood of $S$ is defined by $N(S)=\bigcup_{v\in S} N(v)$, and the closed neighbourhood of $S$ by
$N[S]=N(S)\cup S$. 
 For a set $S \subseteq V$, the subgraph induced by $S$ is denoted by $G[S]$ while the graph $G-S$ is the graph obtained from $G$ by deleting the vertices
in $S$ and all edges incident with $S$. We denote the  minimum and maximum degree among the vertices of $G$  by $\delta (G)$ and $\Delta (G)$, respectively,   or simply by $\delta$ and $\Delta$  if the graph $G$ is clear
from the context. The \textit{girth}  $g(G)$ of a graph $G$ is the length of a shortest cycle 
in $G$. As usual we denote the complement of graph $G$ by $\overline{G}$. 
A set $S$ of vertices in a graph $G$ is a \textit{dominating set}  of $G$ if every vertex of $V(G)\setminus S$ is adjacent to some vertex in $S$, that is $N[S]=V$. The minimum
cardinality of a dominating set of $G$ is the \textit{domination number} of $G$ and  denoted by $\gamma (G)$ or simply $\gamma$.   A $\gamma$-set of $G$ is a set $S$ which  is a dominating set with cardinality $\gamma$.   Here, we state some known results on the domination number which are needed in the next section:
\begin{theorem}{\rm \cite{Haynes1}}\label{domination}
\begin{enumerate}
\item[(i)] For any graph $G$ of order $n$,  $\left\lceil \dfrac{n}{1+\Delta}  \right\rceil \leq \gamma (G) \leq n-\Delta$.
\item[(ii)] If a graph $G$ has $\delta(G) \geq  2$ and $g(G) \geq 7$, then $\gamma (G) \geq \Delta$. \label{partdom}
\item[(iii)] If graph $G$ is disconnected, then $\gamma (\overline{G})\leq 2$.
\end{enumerate}
\end{theorem}

\medskip
A labeling of a simple graph $G$, $\phi : V \rightarrow \{1, 2, \ldots , r\}$, is said to be \textit{$r$-distinguishing},  if no non-trivial  automorphism of $G$ preserves all of the vertex labels. The point of the labels on the vertices is to destroy the symmetries of the
graph, that is, to make the automorphism group of the labeled graph trivial.
Formally, $\phi$ is $r$-distinguishing if for every non-identity $\sigma \in {\rm Aut}(G)$, there exists $x$ in $V$ such that $\phi(x) \neq \phi(\sigma(x))$. The \textit{distinguishing number} of a graph $G$ is defined  by
\begin{equation*}
D(G) = {\rm min}\{r \vert ~ G ~\text{{\rm has a labeling that is $r$-distinguishing}}\}.
\end{equation*} 

This number has defined in  \cite{Albert}. Similar to this definition,  the \textit{distinguishing  index} $D'(G)$ of a graph $G$ has defined as  the least integer $d$ such that $G$ has an
edge labeling with $d$ labels that is preserved only by the identity automorphism of $G$.  If a graph has no nontrivial automorphisms, its distinguishing number is  $1$. In other words, $D(G) = 1$ for the asymmetric graphs.
 The other extreme, $D(G) = \vert V(G) \vert$, occurs if and only if $G = K_n$. The distinguishing index of some examples of graphs was exhibited in \cite{Kali1}. For 
 instance, $D(P_n) = D'(P_n)=2$ for every $n\geq 3$, and 
 $D(C_n) = D'(C_n)=3$ for $n =3,4,5$,  $D(C_n) = D'(C_n)=2$ for $n \geq 6$. It is easy to see that the value $|D(G)-D'(G)|$ can be large. For example $D'(K_{p,p})=2$ and $D(K_{p,p})=p+1$, for $p\geq 4$.  In the sequel, we need the following results:
\begin{theorem}{\rm \cite{K.L. Collins A.N. Trenk,klavzaretal}}\label{disnumbound}  If $G$ is a connected graph with maximum degree $\Delta$, then $D(G) \leq \Delta+1$. Furthermore,
equality holds if and only if $G$ is a $K_n$, $K_{n,n}$, $C_3$, $C_4$ or $C_5$.
 \end{theorem}
 \begin{theorem}{\rm \cite{Kali1}}\label{distindex1}
 If $G$ is a connected graph of order $n \geq 3$, then
$D'(G) \leq \Delta(G)$, unless $G$ is $C_3$, $C_4$ or $C_5$.
\end{theorem}
 \begin{theorem}{\rm \cite{pilsniaknord}}\label{distindex2}
 Let $G$ be a connected graph that is neither a symmetric nor an asymmetric tree. If the maximum degree of $G$ is at least $3$, then $D'(G) \leq \Delta(G) - 1$ unless $G$ is $K_4$ or $K_{3,3}$.
  \end{theorem}

\medskip 
In the next section, we investigate the relationship between the distinguishing number (index) and the domination number of a graph $G$. For any two 
natural numbers  $\gamma$ and $ d$, we present  a connected graph $G$ such that $\gamma(G)=\gamma$ and $D(G)=d$. 
  In Section 3, we propose a problem and state a result about graph $G$ with $D'(G)\leq \gamma(G)$.

\section{$\gamma(G)$  versus $D(G)$ and $D'(G)$}

We begin this section by an observation which is an immediate  consequence of Theorems   \ref{domination} and \ref{disnumbound}.  
\begin{observation}\label{obser1}
Let $G$ be a simple graph of order $n$. We have,
\begin{enumerate}
\item[(i)]  $D(G)- \gamma (G) \leq n-\left\lceil \dfrac{n}{1+\Delta}  \right\rceil$.
\item[(ii)]  $D(G)- \gamma (G) \leq \Delta +1 -\left\lceil \dfrac{n}{1+\Delta}  \right\rceil$.
\item[(iii)] $D(G)\leq n- \gamma (G)$, except for the complete graphs.\label{part3ofobserva}
\item[(iv)] If  $\delta(G) \geq  2$ and $g(G) \geq 7$, then $D(G)\leq \gamma (G)$. 
\item[(v)] If graph $G$ is a connected graph such that $\overline{G}$ is disconnected, then $\gamma (G)\leq D(G)+1$.\label{part5ofobserva}
\item[(vi)]  If graph $G$ is a connected graph of order $n$, then $\dfrac{1}{n-\Delta}\leq \dfrac{D(G)}{\gamma (G)}\leq \dfrac{(\Delta +1)^2}{n}$.
\end{enumerate}
\end{observation}

 \begin{remark}
By Theorems \ref{distindex1} and \ref{distindex2}, we can see that the  Observation \ref{obser1} is true  for the distinguishing index of $G$, $D'(G)$. Note that  
in this case the Part (iii) of Observation \ref{obser1} is true  for all connected graphs,  except $K_3$. 
\end{remark} 

\begin{figure}
	\begin{center}
		\includegraphics[width=0.55\textwidth]{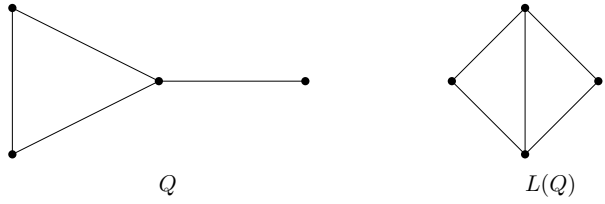}
		\caption{\label{line1} graphs $Q$ and $L(Q)$ of Theorem \ref{autlinegraph}.}
	\end{center}
\end{figure}

 For a simple graph $G$,  the line graph  $L(G)$ is  the graph whose vertices are edges of $G$ and  two edges $e, e' \in   V (L(G)) = E(G)$ are adjacent if they share an endpoint in common.  See Figure \ref{line1}. To study the distinguishing number and the distinguishing  index of  $L(G)$, we need more information about the automorphism group of $L(G)$.
 Let  $\gamma_G : {\rm Aut} (G) \rightarrow {\rm Aut} (L(G))$ be given by $(\gamma_G \phi)(\{u, v\}) = \{\phi(u), \phi(v)\}$ for every $\{u, v\} \in  E(G)$. In \cite{Sabidussi}, Sabidussi proved the following Theorem which we need it later.
  \begin{theorem}{\rm \cite{Sabidussi}}\label{autlinegraph}
 Suppose that $G$ is a connected graph that is not $P_2, Q$, or $L(Q)$ (see Figure \ref{line1}). Then 
$G$ is a group isomorphism, and so ${\rm Aut}(G) \cong {\rm Aut}(L(G))$.
  \end{theorem}
 \begin{theorem}\label{disnumline}
Suppose that $G$ is a connected graph that is not $P_2$ and $L(Q)$ (Figure \ref{line1}). Then $D(L(G))= D'(G)$.
\end{theorem}
\proof  If $G = Q$, then it is easy to see that $D'(Q) = D(L(Q))=2$.  If $G \neq Q$,  First we show that $D(L(G))\leq D'(G)$. Let $c:E(G)\rightarrow \{1, \ldots , D'(G)\}$ be an edge distinguishing labeling of $G$. We define $c': V(L(G)) \rightarrow \{1, \ldots , D'(G)\}$ such that $c'(e) = c(e)$, where $e\in V(L(G))  = E(G)$. The vertex labeling $c'$ is a distinguishing vertex labeling of $L(G)$, because if $f$ is an automorphism of $L(G)$ preserving the labeling, then $c'(f(e))= c'(e)$, and hence $c(f(e))= c(e)$ for any $e\in E(G)$. On the other hand,  by Theorem \ref{autlinegraph}, $f = \gamma_G \phi$ for some automorphism $\phi$ of $G$. Thus from  $c(f(e))= c(e)$ for any $e\in E(G)$, we can conclude that $c(\gamma_G \phi (e))= c(e)$ and so $c(\{\phi(u), \phi(v)\})= c(\{u,v\})$ for every $\{u,v\}\in E(G)$. This means that $\phi$ is an automorphism of $G$ preserving the labeling $c$, and so $\phi$ is the identity automorphism of $G$. Therefore $f$ is the identity automorphism of $L(G)$, and hence $D(L(G))\leq D'(G)$.  For the converse, suppose that  $c:V(L(G))\rightarrow \{1, \ldots , D(L(G))\}$ is a vertex distinguishing labeling of $L(G$). We define $c': E(G) \rightarrow \{1, \ldots ,D(L(G))\}$ such that $c'(e) = c(e)$ where $e\in    E(G)$. The edge labeling $c'$ is a distinguishing edge labeling of $G$. Because if $f$ is an automorphism of $G$ preserving the labeling, then $c'(f(e))= c'(e)$, and hence $c(f(e))= c(e)$ for any $e\in E(G)$. Then, there exist the automorphism $\gamma_G f$ of $L(G)$ such that  $\gamma_G f (\{u,v\})= \{f(u), f(v)\}$ for every $\{u,v\}\in E(G)$, by Theorem \ref{autlinegraph}. Thus from  $c(f(e))= c(e)$ for any $e\in E(G)$, we can conclude that  $c(\{u,v\})= c(\{f(u), f(v)\})= c (\gamma_G f (\{u,v\}))$ for every $\{u,v\}\in E(G)$, which means that  $\gamma_G f$ preserves the distinguishing  vertex labeling of $L(G)$, and hence $\gamma_G f$ is the identity automorphism of $L(G)$. Therefore $f$ is the identity automorphism of $G$, and so $D'(G) \leq D(L(G))$. \qed

By Observation \ref{obser1} (iii) and Theorem \ref{disnumline}, we can obtain a new upper bound for the distinguishing index of a graph using the domination number of its line graph.
 \begin{corollary}
Suppose that $G$ is a connected graph of order at least three and size $m$. Then $D'(G)\leq m-\gamma (L(G))$, except for star graphs.
\end{corollary}

\begin{figure}
	\begin{center}
		\includegraphics[width=1\textwidth]{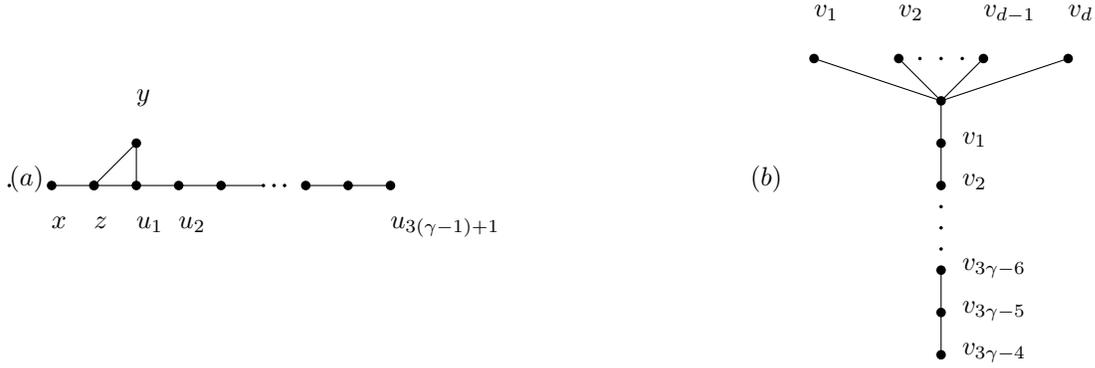}
		\caption{\label{figthm} Graphs  in Theorem \ref{priv}.}
	\end{center}
\end{figure}

Now, we state and prove one of the main result of this paper.  
\begin{theorem}\label{priv}
\begin{enumerate}
\item[(i)] For any two natural numbers $\gamma$ and $ d$, there exists a connected graph $G$ such that $\gamma(G)=\gamma$ and $D(G)=d$.
\item[(ii)] For any two natural numbers $\gamma$ and $ d$, there exists a connected graph $G$ such that $\gamma(G)=\gamma$ and $D'(G)=d$.
\end{enumerate}
\end{theorem}
\proof (i) We consider the following cases:
\begin{itemize}
\item If $2\leq \gamma \leq d$, then we consider the graph $G$ as a tree with central vertex $x$ of degree $\gamma$, and $N(x)=\{x_1,\ldots , x_{\gamma}\}$. Suppose that each vertex $x_i$ is adjacent to $d$ pendant vertices. Since $\gamma \leq d$, so $D(G)=d$. Also, it is clear that the set $\Gamma=\{x_1,\ldots , x_{\gamma}\}$ is the minimum dominating set of $G$.
\item If $\gamma =1$ and $d\geq 2$, then it is sufficient to consider the graph $G$ as the star graph $K_{1,d}$. Hence, $D(G)=d$ and $\gamma (G) =1$.
\item If $\gamma \geq 2$ and $d=1$, then we consider the graph $G$ as the asymmetric graph shown in Figure \ref{figthm} (a). 
\item If $2\leq d \leq \gamma$, then  we consider the graph $G$ as a star graph $K_{1,d}$ with a path of length $3\gamma -3$ attached to the central vertex of  $K_{1,d}$, see Figure \ref{figthm} (b).
\item If $\gamma = d =1$, then we consider the graph $G$ as the join of an asymmetric graph $H$ and $K_1$ such that the graph $H$ does not have  any vertex of degree $|V(H)|-1$. In this case, the graph  $G$ is an asymmetric graph with a vertex of degree $|V(G)|-1$. Hence $D(G)=1$ and $\gamma(G)=1$.
\end{itemize}

The proof of Part (ii) is exactly the same as Part (i).
\qed

The graph $G-v$ is a graph that is made  by deleting the vertex $v$ and all edges connected to $v$ from the graph $G$. The following theorem examine the effects on $D(G)$ and
$D'(G)$ when $G$ is modified by deleting a  vertex of $G$. 
\begin{theorem}{\rm \cite{alikhani}}\label{thmG-v}
Let $G$ be a connected graph of order $n\geq 3$ and $v\in V(G)$. Then we have
\begin{itemize}
\item[(i)] $D(G) -1\leq D(G-v)\leq 2D(G)$.
\item[(ii)] $D'(G) -1\leq D'(G-v)\leq 2D'(G)$.
\end{itemize}
\end{theorem}

By Theorem \ref{thmG-v}, we can obtain the two following bounds for the distinguishing number and the distinguishing  index of a graph based on its domination number.
\begin{corollary}\label{corG-v}
Let $S$ be a $\gamma$-set of  a connected graph $G$ of order $n\geq 3$. We have      
\begin{itemize}
\item[(i)] $D(G) \leq D(G-S)+ \gamma (G)$. 
\item[(ii)]  If the induced subgraph $G-S$ does not have  $K_2$ as its connected component, then $D'(G) \leq D'(G-S)+ \gamma (G)$.
\end{itemize}
\end{corollary}

Now by  Observation \ref{obser1} (iii) and Corollary \ref{corG-v}, we can  prove the following result.
\begin{corollary}\label{corG-v2}
Let $S$ be a $\gamma$-set of  a connected graph $G$ of order $n\geq 3$. We have
\begin{itemize}
\item[(i)] $D(G) \leq \frac{1}{2}(n+D(G-S))$.
\item[(ii)] If the induced subgraph $G-S$ does not have  $K_2$ as its connected component, then  $D'(G) \leq \frac{1}{2}(n+D'(G-S))$. 
\end{itemize}
\end{corollary}

 Pil\'sniak in \cite{pilsniaknord}, showed that the distinguishing index of connected graphs without the graph $K_{1,3}$ as its induced subgraph, i.e., \textit{claw-free graphs}, has the distinguishing index at most three. She proved that the distinguishing index of graphs with a Hamiltonian path, i.e., \textit{tracable graphs}, of order at least 7, is at most two. By this argument and   Corollary \ref{corG-v} (ii), we can state the following result.
 \begin{theorem}
Let $S$ be a $\gamma$-set of  a connected graph $G$ of order $n\geq 3$. We have  
 \begin{enumerate}
 \item[(i)] If $G-S$ is a connected claw-free graph, then $D'(G)\leq \gamma(G)+3$.
 \item[(ii)] If $G-S$ is a connected tracable graph of order at least $7$, then $D'(G)\leq \gamma(G)+2$.
 \end{enumerate}
 \end{theorem}

We close this section by the following result.
\begin{theorem}
Let $S$ be a $\gamma$-set of  a connected graph $G$ of order $n\geq 3$ such that $G-S$ is a connected graph. We have  
\begin{itemize}
\item[(i)] $D(G) \leq n- \gamma(G-S)$.
\item[(ii)]  $D'(G) \leq n- \gamma(G-S)$.
\item[(iii)] $D(G) \leq n -\frac{1}{2}(\gamma(G)+\gamma(G-S))$.
\item[(iv)] If the induced subgraph $G-S$ does not have  $K_2$ as its connected component, then
 $D'(G) \leq n -\frac{1}{2}(\gamma(G)+\gamma(G-S))$ 
\end{itemize}
\end{theorem}
\proof
\begin{itemize}
	\item[(i)]
  By Corollary \ref{corG-v}, we have $D(G) \leq D(G-S)+ \gamma (G)$. On the other hand,  $D(G-S) \leq (n-\gamma (G))-\gamma(G-S)$,  by  Observation \ref{obser1} (iii).
 Hence the result follows.
 \item[(ii)]   The proof is similar to the  proof of Part (i).

\item[(iii)] By using Corollary \ref{corG-v2}, and  by  Observation \ref{obser1} (iii) for graph $G-S$, we have the result.

\item[(iv)]  The proof is similar to the proof of Part (iii).\qed
\end{itemize}

\section{Graphs $G$ with $D'(G)\leq \gamma(G)$}

As we have seen in Observation \ref{obser1}, the distinguishing index of a graph $G$ can be more or less than its domination number. 
With this motivation we propose the following problem. 

\begin{problem}\label{P1}
	\begin{enumerate}
		\item [(i)] Characterize  the connected graphs $G$ with $D(G)=\gamma (G)$ and $D'(G)=\gamma (G)$.
		\item[(ii)] Characterize  the connected graphs $G$ and $H$ with $D'(G)\leq \gamma (G)$ and $\gamma (H) \leq D'(H)$.
		
		\end{enumerate}

\end{problem}

We tried to solve this problem. To state our  result, we need the following theorem which is easy to obtain.   
 \begin{theorem}\label{girth}
 	Let $G$ be a simple graph that is not a forest and has girth at least $5$. Then the complement of $G$ is Hamiltonian.
 \end{theorem}
 
 By  Observation \ref{obser1}$(v)$, we can suppose that both of $G$ and $\overline{G}$ are connected graphs. 

\begin{theorem} \label{char}
	Let  $G$ and $\overline{G}$  are connected graphs.
\begin{enumerate} 
	\item [(i)] If $\gamma (G) \geq \Delta$, then $D'(G) \leq \gamma (G)$, except for graphs $C_3$, $C_4$ and $C_5$. 
	\item[(ii)] If $\gamma (G) <\Delta$    and   $g\geq 5$, then  $D'(\overline{G})\leq 2$.
	\item[(iii)]   If  $\gamma (G) <\Delta$    and     $g=4$, then $D'(\overline{G})\leq 3$.
	\item[(iv)]  If  $\gamma (G) <\Delta$,  $g\leq 3$,  and $\gamma (G) < \Delta \leq n-\Delta$, then $D'(\overline{G})\leq 2$. 
	\item[(v)]  If  $\gamma (G) <\Delta$, $g=3$   and $\delta =1$, then $\gamma(\overline{G})=2$. 
	\item[(vi)]   If $G$ is a tree,  then  $D'(\overline{G})\leq 3$. 
	\end{enumerate}
\end{theorem}
 \proof   
  \begin{enumerate}
 	\item [(i)] It follows from  Theorem \ref{distindex1}.
  	\item[(ii)]  Since  $\overline{G}$ is Hamiltonian, by Theorem \ref{girth} we have  $D'(\overline{G})\leq 2$. 
 	\item[(iii)]  In this case,  the graph $G$ is triangle free, and so  $\overline{G}$ is claw-free. Therefore  $D'(\overline{G})\leq 3$.
 	\item[(iv)]  We have $\delta (\overline{G})\geq n/2$, and so  $\overline{G}$ is Hamiltonian. Therefore  $D'(\overline{G})\leq 2$.
 	\item[(v)] In this case, $\Delta(\overline{G})= n-2$, and so $\gamma(\overline{G})=2$. 
 	\item[(vi)]    Since the graph $G$ is triangle free, the graph  $\overline{G}$ is claw-free, and so $D'(\overline{G})\leq 3$.\qed
 \end{enumerate}

We end this paper with the following remark. 
\begin{remark}  	 	 
Theorem \ref{char} is an elementary answer to Problem \ref{P1}. The exact characterization for this problem remain as an open problem. Note that in Theorem \ref{char}, we have no result for the case  $g=3$ and $\delta \geq 2$. So this case also is open.
 \end{remark}

\end{document}